\documentclass{gtart_h}

\def\ifplaintex{\expandafter\ifx\csname documentclass\endcsname\relax}

\def\gtp{{\mathsurround=0pt\it $\cal G\mskip-2mu$eometry \&\ 
$\cal T\!\!$opology $\cal P\!$ublications}}  

\def\recd{{\small Received:\qua\receiveddate\ifx\reviseddate\relax
\else\qquad Revised:\qua\reviseddate\fi\par}} 

\def\lognumber#1{\def\thelognumber{#1}}
\def\volumenumber#1{\def\thevolumenumber{#1}}
\def\volumeyear#1{\def\thevolumeyear{#1}}
\def\papernumber#1{\def\thepapernumber{#1}}
\def\pagenumbers#1#2{\def\startpage{#1}\def\finishpage{#2}}
\def\published#1{\def\publishdate{#1}}

\def\received#1{\def\receiveddate{#1}}

\def\accepted#1{\def\accepteddate{#1}}

\def\asciiemail#1{\def\theasciiemail{#1}}

\long\def\asciiabstract#1{\long\def\theasciiabstract{#1}}

\let\\\par\let\thelognumber\relax\let\thevolumenumber\relax
\let\thepapernumber\relax\let\thevolumeyear\relax\let\startpage\relax
\let\finishpage\relax\let\publishdate\relax\let\receiveddate\relax
\let\reviseddate\relax\let\accepteddate\relax\let\theasciititle\relax
\let\theasciiauthors\relax
\let\theasciiabstract\relax

\let\theasciiemail\relax

\ifplaintex
\font\logobig=cmssbx10 scaled 3836
\font\logomed=cmssbx10 scaled 2557
\else
\font\logobig=cmssbx10 scaled 4200
\font\logomed=cmssbx10 scaled 2800
\fi

\long\def\makeagttitle{   
\count0=\startpage
\agt\hfill      
\hbox to 45truept{\vbox to 0pt{\vglue -13truept{\logomed A\kern -.37em{\logobig 
T}\kern -.38em G}\vss}\hss}
\break
{\small Volume \thevolumenumber\ (\thevolumeyear)
\startpage--\finishpage\nl
Published: \publishdate}

\vglue .25truein

{\parskip=0pt\leftskip 0pt plus
1fil\def\\{\par\smallskip}{\Large\bf\thetitle}\par\medskip} \vglue
0.05truein

{\parskip=0pt\leftskip 0pt plus 1fil\def\\{\par}{\sc\theauthors}
\par\medskip}%
 
\vglue 0.03truein 

{\small\leftskip 25truept\rightskip 25truept{\bf Abstract}\stdspace\theabstract

{\bf AMS Classification}\stdspace\theprimaryclass
\ifx\thesecondaryclass\relax\else; \thesecondaryclass\fi\par
{\bf Keywords}\stdspace \thekeywords\par}\vglue 7truept

}   

\ifplaintex
\font\phead=cmsl9 scaled 950
\font\pnum=cmbx10 scaled 913
\font\pfoot=cmsl9 scaled 950
\headline{\vbox to 0pt{\vskip -4.5mm\line{\small\phead\ifnum
\count0=\startpage ISSN 1472-2739 (on-line) 1472-2747 (printed)
\hfill {\pnum\folio}\else\ifodd\count0\def\\{ }%
\ifx\theshorttitle\relax\thetitle\else\theshorttitle\fi\hfill{\pnum\folio}
\else\def\\{ and }{\pnum\folio}\hfill\ifx\theshortauthors\relax\theauthors
\else\theshortauthors\fi\fi\fi}\vss}}
\footline{\vbox to 0pt{\vglue 0mm\line{\small\pfoot\ifnum\count0=\startpage
\copyright\ \gtp\hfill\else
\agt, Volume \thevolumenumber\ (\thevolumeyear)\hfill\fi}\vss}}
\else
\font=cmsl9 scaled 1050
\font\lnum=cmbx10 
\font=cmsl9 scaled 1050
\makeatletter
\def\@oddhead{{\small\lhead\ifnum\count0=\startpage ISSN 1472-2739 
(on-line) 1472-2747 (printed)\hfill {\lnum\number\count0}\else\ifodd\count0
\def\\{ }\ifx\theshorttitle\relax \thetitle \else\theshorttitle\fi\hfill
{\lnum\number\count0}\else\def\\{ and }{\lnum\number\count0}
\hfill\ifx\theshortauthors\relax 
\theauthors\else\theshortauthors\fi\fi\fi}}\def\@evenhead{\@oddhead}
\def\@oddfoot{\small\lfoot\ifnum\count0=\startpage\copyright\ \gtp\hfill\else
\agt, Volume \thevolumenumber\ (\thevolumeyear)\hfill\fi}
\def\@evenfoot{\@oddfoot}
\makeatother
\fi
\let\maketitlepage\makeagttitle
\let\makeshorttitle\maketitlepage
\let\maketitle\maketitlepage

\newwrite\gtoutfile
\long\gdef\makeheadfile{  
{\def\\{, }\def\s{ }
\immediate\openout\gtoutfile head.xxx
\immediate\write\gtoutfile{Proxy-for: \ifx\theasciiauthors\relax
\theauthors\else\theasciiauthors\fi\s<\ifx\theasciiemail\relax\theemail\else\theasciiemail\fi>}
\immediate\write\gtoutfile{\noexpand\\}
\immediate\write\gtoutfile{Authors: \ifx\theasciiauthors\relax
\theauthors\else\theasciiauthors\fi}
{\def\\{ }\immediate\write\gtoutfile{Title: \ifx\theasciititle\relax
\thetitle\else\theasciititle\fi}}
\immediate\write\gtoutfile{Subj-class: GT or SG, GR etc}
\immediate\write\gtoutfile{MSC-class: \theprimaryclass\ifx\thesecondaryclass\relax\else, \thesecondaryclass\fi}
\immediate\write\gtoutfile{Journal-ref: Algebr. Geom. Topol. \thevolumenumber\s
(\thevolumeyear) \startpage-\finishpage}
\immediate\write\gtoutfile{Comments: Published by Algebraic and
Geometric Topology at}
\immediate\write\gtoutfile{\s\s\s  http://www.maths.warwick.ac.uk/agt/AGTVol\thevolumenumber/agt-\thevolumenumber-\thepapernumber.abs.html}
\immediate\write\gtoutfile{\noexpand\\}
\immediate\write\gtoutfile{}
\ifx\theasciiabstract\relax
\immediate\write\gtoutfile{\theabstract}\else
\immediate\write\gtoutfile{\theasciiabstract}\fi
\immediate\write\gtoutfile{}
\immediate\write\gtoutfile{\noexpand\\}
\immediate\write\gtoutfile{}
\immediate\closeout\gtoutfile}}  

\def\maketitlepage{\makeagttitle\makeheadfile}
\let\makeshorttitle\maketitlepage
\let\maketitle\maketitlepage

\lognumber{54}
\volumenumber{4}
\volumeyear{2004}
\papernumber{54}
\pagenumbers{1253}{1272}
\received{28 June 2004} 
\accepted{7 December 2004}
\published{24 December 2004}

\usepackage{amssymb, amsmath, latexsym, amscd}
\input epsf
\def\relabelbox{%
  \hbox\bgroup%
}%
\def\endrelabelbox{%
\egroup%
}%
\def\relabel #1#2 {%
  \special{ps:/a {} def}%
  \smash{\rlap{#2}}%
}%
\def\adjustrelabel <#1,#2> #3#4 {%
  \special{ps:/a {} def}%
  \smash{\rlap{\kern #1 \raise #2\hbox{#4}}}%
}%
\def\extralabel <#1,#2> #3 {\smash{\rlap{\kern #1 \raise #2\hbox{#3}}}}%

\let\Bbb\mathbb   
\newtheorem{theorem}{Theorem}
\newtheorem{corollary}{Corollary}
\newtheorem{lemma}{Lemma}
\newtheorem{proposition}{Proposition}

\def\ontop#1{\buildrel #1\over \to}

\def\G{\Gamma}

\newcommand{\bdry}{\partial} 
\newcommand{\Diff}{{\it Diff}}
\def\Z{\Bbb Z}

\begin{document} 

\title{Homology stability for outer automorphism\\groups of free groups}\author{Allen Hatcher\\Karen Vogtmann}
\address{Department of Mathematics, Cornell 
University\\Ithaca, NY 14853-4201, USA}
\gtemail{\mailto{hatcher@math.cornell.edu}{\rm\qua 
and\qua}\mailto{vogtmann@math.cornell.edu}}                
\asciiemail{hatcher@math.cornell.edu, vogtmann@math.cornell.edu} 

\begin{abstract} 
We prove that the quotient map from $Aut(F_n)$ to $Out(F_n)$ induces
an isomorphism on homology in dimension $i$ for $n$ at least
$2i+4$. This corrects an earlier proof by the first author and
significantly improves the stability range. In the course of the
proof, we also prove homology stability for a sequence of groups which
are natural analogs of mapping class groups of surfaces with
punctures. In particular, this leads to a slight improvement on the
known stability range for $Aut(F_n)$, showing that its $i$th homology
is independent of $n$ for $n $ at least $2i+2$.  
\end{abstract}
 
\asciiabstract{%
We prove that the quotient map from Aut(F_n) to Out(F_n) induces an
isomorphism on homology in dimension i for n at least 2i+4. This
corrects an earlier proof by the first author and significantly
improves the stability range. In the course of the proof, we also
prove homology stability for a sequence of groups which are natural
analogs of mapping class groups of surfaces with punctures. In
particular, this leads to a slight improvement on the known stability
range for Aut(F_n), showing that its i-th homology is independent of n
for n at least 2i+2.}


\maketitle

\section {Introduction}

In this paper we present a proof that the quotient map $ Aut(F_n) \to Out(F_n) $ 
induces an isomorphism on the integral homology groups $H_i $ for $n \ge 2i+ 4$. 
 We also show that the natural map
$ Aut(F_n)
\to Aut(F_{n+1})$ induces an isomorphism on $H_i$ for $ n \ge 2i+2 $, a slight 
improvement on  the result in  \cite{HV}.  It follows that
 $H_i(Out(F_n)) $ is independent of $n$ for $n \ge 2i+4$. 

These stable homology groups have a significant amount of non-trivial torsion
since they contain the homology of $\Omega^\infty S^\infty$ as a direct summand,
 as observed in \cite{H}.  In contrast,  it is conjectured that they are rationally trivial,
though a few non-trivial unstable rational classes have been found
 \cite{CoVo, Gerlits, HV}.   One of these non-trivial classes, in $H_7(Aut(F_5);{\bf
Q})$,
 becomes trivial under the map to $H_7(Out(F_5);{\bf Q})$, showing that this map is
not always an
 isomorphism on homology \cite{Gerlits}.

The first proof of homology stability for $Aut(F_n)$ and $Out(F_n)$, with a much
weaker dimension range than we obtain here, appeared in the first author's 1995 paper
\cite{H}.  However, in 2003 a mistake was discovered in the
spectral sequence arguments at the end of that paper (in the last
sentence on p.56, where the diffeomorphism in question may permute boundary spheres and
hence not induce an inner automorphism as claimed). In the meantime an independent proof for the
isomorphisms
$H_i(Aut(F_n)) \cong H_i(Aut(F_{n+1}))$, with a much better dimension range, 
had appeared in 1998 in \cite{HV}. This left the other isomorphisms $H_i(Out(F_n))
\cong H_i(Aut(F_n))$ unproved, and the present paper fills that gap, at the same time
giving a better dimension range than the one originally claimed.

The groups $Out(F_n)$ and $Aut(F_n)$ are the cases $s =0,1$ of a sequence of groups
$\Gamma_{n,s}$ which are analogs of mapping class groups of surfaces with $s$ 
punctures.  We will show  that $ H_i(\Gamma_{n,s}) $ is
independent of $n$ and $s$ for $n \ge 2i+4$. A more general family of groups
$ \Gamma^k_{n,s}$ analogous to mapping class groups of surfaces with $s$ punctures and
$k$ boundary components is studied in \cite{HW}, where homological stability for $
\Gamma^k_{n,s}$ is proved in the weaker range $ n \ge 3i+3$. In particular this gives a
different proof of stability for both $Aut(F_n)$ and $Out(F_n)$ in this weaker range.
The groups $
\Gamma^k_{n,1} $ are the groups $A_{n,k}$ of \cite{JW}, and they are used in \cite{Wahl} to
study the relation between mapping class groups of surfaces and $Out(F_n)$. 

 The present
paper relies heavily on the geometric results of \cite{H} concerning the connectivity of various
sphere complexes, and on the Degree Theorem of \cite{HV}. The spectral sequence
arguments given here are new; they completely replace those in
\cite{H} and can also be used as an alternative to the spectral sequence arguments in \cite{HV}.

\rk{Acknowledgment}  The second author was partially supported by NSF grant DMS 0204185.

\section{Tools and Outline of proof}

Let $M_{n,s}$ be the compact $3$--manifold obtained from the connected sum of $n$ copies
of
$S^1\times S^2$ by deleting the interiors of $s$ disjoint balls. When $ n=0 $
we take $M_{0,s}$ to be $S^3$ with the interiors of $s$ disjoint balls removed.  The
$s$ boundary spheres of $M_{n,s}$ are denoted $ \bdry_0, \cdots, \bdry_{s-1} $,
considering
$
\bdry_0$ as the base sphere. Let $\Diff(M_{n,s})$ be the group of
orientation-preserving diffeomorphisms of
$M_{n,s}$ which fix the boundary pointwise. The mapping class group
$\pi_0\Diff(M_{n,s})$ has a normal subgroup generated by Dehn twist diffeomorphisms
along embedded 2--spheres, and the
 quotient of $\pi_0 \Diff(M_{n,s})$ by this subgroup will be denoted $\Gamma_{n,s}$.
A theorem of Laudenbach \cite{Laudenbach} states that $\Gamma_{n,0}$ is isomorphic to
$Out(F_n)$, and $\Gamma_{n,1}$ to $Aut(F_n)$.

An inclusion $ M_{n,s} \hookrightarrow M_{m,t} $ induces a homomorphism $
\Gamma_{n,s} \to \Gamma_{m,t} $ by extending diffeomorphisms by the identity map
outside $ M_{n,s}$.  We will be interested in three special cases:
\begin{enumerate}
\item {\it Adding a punctured handle}.   The map $\alpha\co  \Gamma_{n,s}  \to \Gamma_{n+1,s}$,
$ s \ge 1 $, is induced by gluing a copy of $M_{1,2}$ to the base sphere $\bdry_0$; the
puncture in
$M_{1,2}$ which is not used in the gluing is the new base sphere. 
\item  {\it Adding a tube}.  The map $\beta\co \Gamma_{n,s} \to \Gamma_{n+1,s-2}$, $ s \ge 2
$,   is induced by gluing the boundary components of a copy of $M_{0,2}$ to the last two
boundary components of
$M_{n,s}$
\item {\it Adding pants}.  
The map $\mu\co \Gamma_{n,s} \to \Gamma_{n,s+1}$, $ s \ge 1$, is induced by gluing one boundary
component of a copy of
$M_{0,3}$ to $\bdry_0$. One of the unglued boundary components of $M_{0,3}$ is
the new base sphere.  \end{enumerate}

\begin{figure}[ht]
\centerline
{\relabelbox\small\epsfxsize 3truein
\epsfbox{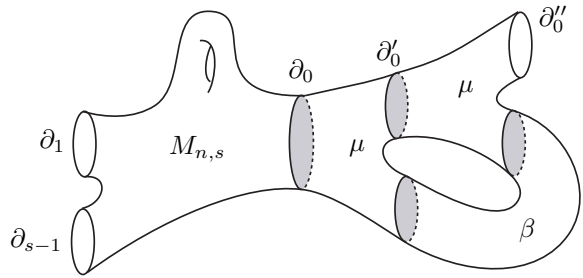}
\relabel{m}{$\mu$}
\relabel{b}{$\beta$}
\relabel{b1}{$\partial_1$}
\relabel{b2}{$\partial_{s-1}$}
\relabel{b3}{$\partial_0^{\prime\prime}$}
\relabel{b5}{$\partial_0^{\prime}$}
\relabel{b4}{$\partial_0$}
\relabel{m1}{$\mu$}
\relabel{Mns}{$M_{n,s}$}
\endrelabelbox}
\caption{$\alpha=\beta\circ \mu^2$}
\label{mmb}
\end{figure}

We consider the following commutative diagram of maps (see Figure ~\ref{mmb}):
 $$
\begin{CD}
 \G_{n,s} @>{\alpha}>> \G_{n+1,s}\\
@VV{\mu}V @AA{\beta}A \\
\G_{n,s+1} @>{\mu}>> \G_{n,s+2}
\end{CD}
$$
We first show that $\alpha$ is an isomorphism  on $H_i$ for $n\geq 2i+2$.  This
generalizes and improves the theorem in \cite{HV}  that  $Aut(F_{n})\to Aut(F_{n+1})$
is an isomorphism on $H_i$  for $n\geq 2i+3$.   We  prove this by working with the
spectral sequence associated to the action of $\Gamma_{n,s}$ on a simplicial complex
$W_{n,s}$.  A $k$--simplex of $W_{n,s}$ is a collection of $k+1$ spheres whose complement  in $M_{n,s}$ is connected, together with an additional {\it enveloping} sphere which separates these spheres and the base sphere from the rest of the manifold. 

It follows by a formal argument that $\beta$ and $\mu$ induce isomorphisms on  $H_i$
for $n\ge2i+2$ when
$s> 0$. However, we need to show that $ \beta$ induces isomorphisms on homology when
$s=0$ as well, and for this we use a similar spectral sequence argument, replacing the
complex $W_{n,s}$ by a slightly different complex of spheres $X_{n,s}$ that was
considered
 in \cite{H}. This argument requires $ n \ge 2i+3$.
\begin{figure}[ht]
\centerline
{\relabelbox\small\epsfxsize 3truein
\epsfbox{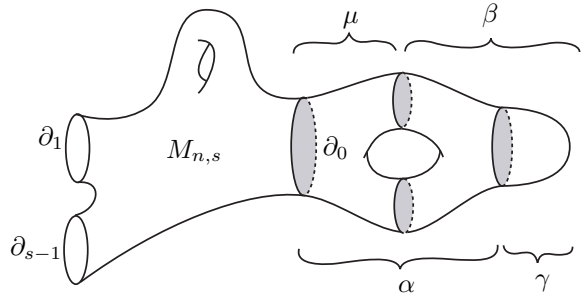}
\relabel{m}{$\mu$}
\relabel{b}{$\beta$}
\relabel{b1}{$\partial_1$}
\relabel{bs}{$\partial_{s-1}$}
\relabel{b0}{$\partial_0$}
\relabel{a}{$\alpha$}
\relabel{g}{$\gamma$}
\relabel{Mns}{$M_{n,s}$}
\endrelabelbox}
\caption{$\gamma\circ\alpha=\beta\circ \mu$}\label{mbag}
\end{figure}

Finally we consider the commutative diagram
$$
\begin{CD}
\Gamma_{n,1}@>{\alpha}>>\Gamma_{n+1,1}\\
@VV{\mu}V @VV{\gamma}V\\
\Gamma_{n,2}@>{\beta}>>\Gamma_{n+1,0}
\end{CD}
$$
where $ \gamma $ is induced by filling in the boundary sphere with a ball (see Figure~\ref{mbag}). 
 Thus $\gamma$ is the
quotient map $Aut(F_{n+1}) \to Out(F_{n+1})$. Since $ \alpha$, $\beta$, and $\mu$
induce isomorphisms on $H_i$ for $n \geq 2i+3$, so does $\gamma$. Replacing $n+1$ by
$n$, this says that the map $Aut(F_n) \to Out(F_n)$ induces an isomorphism on $H_i$
for $ n
\ge 2i+4$.

The main input needed to make the spectral sequences work is the fact that the
complexes
$W_{n,s}$ and $X_{n,s}$ are highly connected, as are their quotients by the action of
$\G_{n,s}$. For $X_{n,s}$ and its quotient this was shown in \cite{H}. The
connectivity of $W_{n,s}$ will follow from an easy generalization of the main
technical result of \cite{HV}, the Degree Theorem. The quotient
$W_{n,s}/\G_{n,s}$ is combinatorially the same as $X_{n,s}/\G_{n,s}$.

\section{Sphere systems and graphs}

\begin{figure}[ht]
\centerline
{\relabelbox\small\epsfxsize 3truein
\epsfbox{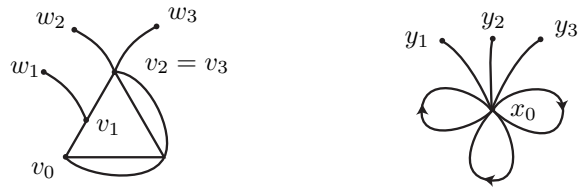}
\relabel{v0}{$v_0$}
\relabel{v1}{$v_1$}
\relabel{v23}{$v_2=v_3$}
\relabel{w1}{$w_1$}
\relabel{w2}{$w_2$}
\relabel{w3}{$w_3$}
\relabel{y1}{$y_1$}
\relabel{y2}{$y_2$}
\relabel{y3}{$y_3$}
\relabel{x0}{$x_0$}
\endrelabelbox}
\caption{A thorned graph with 3 thorns and the standard thorned rose $R_{3,4}$}
\label{thorns}
\end{figure}

We will show that the complex $W_{n,s}$ is  highly connected by comparing it to a complex of
metric graphs with $s$ distinguished points, which are allowed to coincide. The first distinguished point is the basepoint
$v_0$, and the rest will be labeled 
$v_1,\cdots,v_{s-1}$.  The metric is determined by the lengths of the edges, and is normalized so that the sum of these lengths is 1.    It is convenient
also to attach a free edge  of some fixed positive length  to each
$v_i$, $i>0$, and we label the free vertex of this edge by
$w_i$.  
The resulting graph is called a {\it thorned graph}, with thorns the $s-1$ 
free edges attached at $v_1,\cdots,v_{s-1}$ (see Figure ~\ref{thorns}). 

The {\it standard thorned rose} $R_{n,s}$ is a bouquet of $n$ circles with vertex $x_0$
together with
$s-1$ free edges attached at $x_0$, with free vertices labeled $y_1,\cdots, y_{s-1}$.   Each
(oriented) petal of the rose is identified with a generator of the free group $F_n$. There are
natural inclusions $R_{n,s}\subseteq R_{m,t}$ for all $m\geq n$ and $t\geq s$.

By a {\it marking}  of  a thorned graph $(\Gamma, v_0,\cdots, v_{s-1})$ we will mean  a  homotopy
equivalence $g\colon R_{n,s}\to \Gamma$ which sends $x_0$ to $v_0$ and $y_i$ to $w_i$ for
all $i$.    Two marked graphs $(g,\Gamma, v_0,\cdots, v_{s-1})$ and $(g',\Gamma',
v'_0,\cdots, v'_{s-1})$ are {\it equivalent} if there is an isometry
$h\colon\Gamma\to\Gamma'$ sending $v_i$ to $v'_i$ for all $i$ with $g'\simeq h\circ g$
by a homotopy which is constant on the $v_i$.   The lengths of the edges which are not
thorns give barycentric coordinates for an open simplex associated to the marked
thorned graph $(g, \Gamma, v_0,\cdots, v_{s-1})$.  Codimension one faces of this
simplex correspond to marked thorned graphs obtained by collapsing an edge of $\Gamma$
which is not a loop or a thorn.  The space obtained from the disjoint union of all of
the open simplices associated to marked thorned graphs of rank $n$ with $s-1$ thorns by
attaching them together according to the face relations is denoted $A_{n,s}$.

Recall from \cite{H} the complex $S(M)$ of sphere systems in $M = M_{n,s}$. A vertex of this complex
is an isotopy class of embedded spheres in $M$ that do not bound a ball and are not
isotopic to a sphere of $\bdry M$. A set of $k+1$ vertices of $S(M)$ spans a
$k$--simplex if these vertices are represented by disjoint spheres,
no two of which are isotopic. A point of $S(M)$ can be thought of as a weighted system
of spheres, the weights being given by the barycentric coordinates of the point. A
system
$S$ is said to be {\it simple} if each component of $ M-S$ is simply-connected. We let
$ SS(M) $ be the subspace of $S(M)$ formed by weighted simple systems with strictly
positive weights. This is a union of open simplices of $S(M)$.

\begin{figure}[ht]
\centerline
{\relabelbox\small\epsfxsize 3truein
\epsfbox{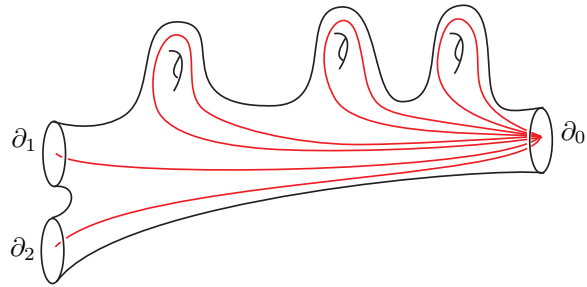}
\relabel{b1}{$\partial_1$}
\relabel{b2}{$\partial_{2}$}
\relabel{b0}{$\partial_0$}
\endrelabelbox}
\caption{$R_{3,3}$ embedded in $M_{3,3}$}\label{embedded}
\end{figure}

As in \cite{H, HV} which dealt with the cases $s\le1$, there is a homeomorphism $ \Phi\colon
SS(M) \to A_{n,s}  $, defined in the following way. Given a sphere
system in $SS(M)$,  thicken each sphere to $S^2\times I$, and also thicken the
boundary spheres other than $\bdry_0$.  Map $M$ to a graph by projecting each
$S^2\times I$ onto the interval
$I$ (an edge of the graph), and each component of the closure of $M$ minus the
thickened spheres to a point (a vertex of the graph). Thus the graph is a quotient
space of $M$.  Each boundary sphere except
$\bdry_0$ will correspond to a thorn of the graph. We assign lengths to the
non-thorn edges according to the weight of the corresponding sphere.  To
get a marking, we fix a particular simple sphere system with connected complement. 
Its associated graph is the standard thorned rose $R_{n,s}$, and this can be embedded
in
$M$ so that the embedding induces an isomorphism on $ \pi_1 $, and so that the
basepoint of the image is located in $\bdry_0$ and the thorns run out to the other
boundary spheres (see Figure~\ref{embedded}). The marking of a sphere system $S$ is
given by this embedding followed by the quotient map from $M$ to the graph associated
to $S$.  

It is shown in the Appendix of \cite{H} that $\Phi$ is a homeomorphism when $ s = 0 $, and the same
argument works for $s>0$ once we give an alternative description of $\Gamma_{n,s}$ in terms of
homotopy equivalences of graphs. Let $ G_{n,s} $ be the group of path-components of the space of
homotopy equivalences of the standard thorned rose $ {R = R_{n,s}} $ that fix the $s$
distinguished points. Such a homotopy equivalence is actually a homotopy equivalence relative to the $s$
distinguished points, by Proposition~0.19 of \cite{HAT}. This shows in particular that
$G_{n,s}$ is a group.

There is a map $ \varphi\colon\Gamma_{n,s} \to G_{n,s} $ sending a
diffeomorphism $ f\colon M \to M $ to the composition  
$$
R \ontop i M \ontop f M \ontop r R 
$$
 where $i$ is the standard inclusion and $r$ is a retraction onto $i(R)$
inducing an isomorphism on $ \pi_1$. Since Dehn twists along $2$--spheres in $M$ have
no effect on the homotopy class of $f\circ i$, the map $\varphi$ is well-defined. To
check that
$\varphi$ is a homomorphism we need to see that a composition
$$
R  \ontop i M \ontop f M \ontop r R  \ontop i M \ontop g M \ontop r R 
$$
can be homotoped so as to eliminate the middle two maps. It is not true that $i\circ r$ is homotopic to the
identity, but  the composition $f\circ i$ can be homotoped, fixing
the distinguished points, to have image in $ i(R)$, so that $i\circ r$ {\it is} the
identity on this image. This is a consequence of the fact that $\pi_1(M,i(R))=0$ since
$i$ induces a surjection on $\pi_1$.

\begin{proposition} The homomorphism $ \varphi \colon \Gamma_{n,s} \to G_{n,s} $ is an
isomorphism.
\end{proposition}
 
\begin{proof}  When $ s \le 1$ this is a theorem of Laudenbach \cite{Laudenbach}. For the general case we do
induction on $s$. For convenience, let us replace $\Diff(M_{n,s})$ by $\Diff_{n,s}$, 
the orientation-preserving diffeomorphisms of $ M_{n,0}$ that fix $s$ distinguished
points, so $\pi_0\Diff_{n,s}$ is $\pi_0\Diff(M_{n,s})$ with Dehn twists
along
 the boundary spheres of $M_{n,s}$ factored out. By evaluating elements of
$\Diff_{n,s-1}$ at the
$s$-th distinguished point we obtain a fibration $ \Diff_{n,s} \to \Diff_{n,s-1} \to
M'$ where $M'$ is the complement of the first $s-1$ distinguished points in $M_{n,0}$.
The last four terms of the first row of the following commutative diagram are part of
the long exact sequence of homotopy groups for this fibration:
$$
\begin{CD}
1 @>>> F_n @>>> \pi_0\Diff_{n,s}@>>>\pi_0\Diff_{n,s-1} @>>> 1 \\
@.@| @VV\psi_{s}V@VV\psi_{s-1}V@. \\
1 @>>> F_n @>>> G_{n,s}@>>>G_{n,s-1} @>>> 1 \\
\end{CD}
$$
The second row is derived in a similar way. We identify $G_{n,s-1}$ with $\pi_0$ of 
the space of homotopy equivalences of $R_{n,s}$ which fix the basepoint and the first
$s-2$ thorn tips, but are allowed to do anything on the last thorn.  Restricting such
a homotopy equivalence to the endpoint of the last thorn gives a map to  $R_{n,s}$,
which is a fibration with fiber consisting of homotopy equivalences fixing all the
distinguished points.  The long exact sequence for this fibration gives the second
row. The initial $1$ in this row is
$\pi_1$ of the total space of the fibration, which is zero
since
$\pi_2$ of a graph is trivial. The   maps $\psi_s$ and $\psi_{s-1}$ are the
compositions of the quotient maps to $\G_{n,s}$ and $\G_{n,s-1}$ with 
$\varphi$. The diagram then commutes, and this yields the initial
$1$ in the first row.

A standard diagram chase shows that the kernels of $\psi_s$ and $\psi_{s-1}$   are
isomorphic. By induction we may assume that the kernel of $\psi_{s-1}$ is the subgroup generated by
twists along
$2$--spheres, so it follows that the same is true for $\psi_s$. Thus if we replace each
$\pi_0\Diff$ by $\G$ we still have a commutative diagram of short exact sequences. The five lemma then
finishes the induction step.
\end{proof}

It is shown in Theorem 2.1 of   \cite{H} that $S(M)$ is contractible if $n>0$,  and
we claim that when this contraction is restricted to the subspace $SS(M) $, it 
gives a contraction of $SS(M)$. Verifying this claim is simply a matter of
following through the proof for $S(M)$ and checking that if one starts in $SS(M)$ then one
never leaves it, using the observation that surgery on a simple system always
produces a simple system. As the proof takes two or three pages, we will not repeat
it here. The proof for $SS(M)$ is actually somewhat simpler since
Lemma~2.2 in the proof for $S(M)$ is not necessary for $SS(M)$, in view of the
fact that surgery takes simple systems to simple systems.

When $n = 0$, the space $S(M)=SS(M)$ is not contractible, but is instead homotopy
equivalent to a wedge of spheres of dimension $s-4$. This fact will not be needed
here, so we will not give a proof.

Since $A_{n,s}$ is homeomorphic to $SS(M) $, it  is also contractible when $ n > 0
$.

\section{The degree theorem with distinguished points}

We will need a
variation of the Degree Theorem of \cite{HV}. In that paper, the degree of a basepointed graph of rank $n$ was
defined to be $2n-|v_0|$, where $|v_0|$ is the valence of the basepoint $v_0$.  The Degree
Theorem states that a $k$--parameter family of graphs in  $A_{n,1}$ can be deformed
into the subspace consisting of  graphs of degree at most $k$ by a deformation during which
degree is non-increasing.  Since $A_{n,1}$ is contractible, this implies that  the subspace
of degree $k$ graphs is $(k-1)$--connected.

We define the  degree of a thorned graph $(\Gamma,v_0,\cdots,v_{s-1})$   to be $2n+s-1-|v_0|$.  The thorns  attached at the basepoint contribute to its valence, so that
 the standard thorned rose has degree 0.  

\begin{lemma} If $k<n/2$, then a degree $k$ thorned graph of rank $n$   has at least one loop
at the basepoint. (Here a loop at the basepoint has no thorns attached to its interior.)
\end{lemma}

\begin{proof} The proof is an Euler characteristic argument identical  to the case of 
graphs with no thorns (see \cite{HV}, Lemma~5.2).   Here is how it goes:  Suppose $\Gamma$ is a
thorned graph of rank $n$ and degree $k$, with basepoint $v_0$.  By definition,
$k=2n+s-1-|v_0|$.   Blowing up if necessary, we may assume that all vertices other than $v_0$
are trivalent or univalent. For such a graph the degree equals the number of
 vertices other than $v_0$ that are trivalent since collapsing an edge from $v_0$ to a
trivalent vertex decreases degree by one, and this can be repeated until one has the
standard thorned rose.
 
 Let $\Gamma_1$ denote the subgraph of $\Gamma$ spanned by all vertices other than the 
basepoint.  If $\Gamma$ has no loops at the basepoint, then $E(\Gamma)=E(\Gamma_1) + |v_0|$.
We calculate the Euler characteristic of $\Gamma$ using this fact and the observation that
$V(\Gamma_1)=k+s-1$:
   \begin{align*}
1-n=\chi(\Gamma)&=V(\Gamma)-E(\Gamma)\\
&=\bigl(1+V(\Gamma_1)\bigr)-\bigl(E(\Gamma_1)+|v_0|\bigr)\\
&= (1+k+s-1) - \bigl(E(\Gamma_1) +(2n+s-1-k)\bigr)
\end{align*}
This simplifies to $ n=2k-E(\Gamma_1)$, hence $n\le 2k$.  This was under the assumption that
$\Gamma$ has no loops at the basepoint, so if
$n>2k$ there must be at least one loop at the basepoint.
\end{proof}

We denote by $A_{n,s}^k$   the subspace of $A_{n,s}$ consisting of thorned graphs of degree
at most
$k$.  

The version of the Degree Theorem that we will need is the following.

\begin{theorem} A map $f_0\colon D^k\to A_{n,s}$ is homotopic to a map $f_1\colon D^k\to
A_{n,s}^k$  by a homotopy during which degree decreases monotonically. 
\end{theorem}

\begin{proof}  The proof in \cite{HV} of the Degree Theorem applies directly to thorned graphs, and we refer the reader there for details.  Here is a sketch of the proof. 

We consider the natural  Morse function on each graph measuring distance to the basepoint.  Small initial segments of the edges below each critical point form the  {\it $\epsilon$--cones} of the graph.  The thorns all point  upwards  with respect to the Morse function, so are never part of any $\epsilon$--cones.  
  
The $k$--parameter family of marked graphs is initially perturbed, by varying the edge-lengths of the graphs, to remove all but $k$ critical points (counted with multiplicity) of each graph in the family.  If the graphs have thorns, the remaining critical points may or may not occur at the base of thorns.  

The deformation then proceeds by two operations,  which may need to be repeated
several times.  First is {\it canonical splitting}, which moves non-critical vertices
down towards the basepoint until they hit a critical vertex.  In this process, the
non-critical bases of thorns move downward just as any other vertex would.  The second
operation is {\it sliding in the $\epsilon$--cones,} which perturbs the attaching
points of edges (including thorns) which come down into a critical vertex.  The
perturbation moves attaching points further down the $\epsilon$--cones, so they can
then be canonically split further down the graph.  After these two processes are
repeated sufficiently often, the degree of each
graph in the family is at most~$k$.
\end{proof}

\section{Graphs with loops at the basepoint}

We will prove that $W_{n,s}$ is highly connected by  comparing it to the 
subspace  $L_{n,s}$ of $A_{n,s}$ consisting of graphs with at  least one loop at the basepoint and
no separating edges other than thorns.
\begin{theorem} The subspace $L_{n,s}$ of $A_{n,s}$  is $(n-3)/2$--connected.  
\end{theorem} 

\begin{proof}
Note first that separating edges which are not thorns can always be eliminated by shrinking
them uniformly to points, so this condition in the definition of $L_{n,s}$ is easy to achieve
and we will say no more about it in the proof.

 Suppose $k<n/2$, and take a map of $S^{k-1}$ into $L_{n,s}$. Since $A_{n,s}$ is contractible,
this can be extended to a map $f_0\colon D^k\to A_{n,s}$.   By the Degree Theorem, this is
homotopic to a map $f_1\colon D^k\to  A_{n,s}^k$ by a homotopy which monotonically decreases
degree.  By the preceding Lemma, every graph of degree at most $n/2$ has at least one loop at
the basepoint, so that $f_1(D^k)\subset L_{n,s}$.
 
 The proof of the Degree Theorem shows that  an edge which
terminates at the basepoint at time $t$ remains at the basepoint for all $t'>t$.   In particular, the number of loops at the
basepoint can only increase. Thus the image of
$S^{k-1}$ stays in $L_{n,s}$ for all times $t$ and gets filled in at time
$t=1$ by a
$D^k$ in $L_{n,s}$, so $L_{n,s}$ is $(k-1)$--connected. 
\end{proof}

Let $SL_{n,s}$ denote the spine of $L_{n,s}$, the simplicial complex associated to the poset
of simplices of $L_{n,s}$. Thus $SL_{n,s}$ is a subcomplex of the spine of 
$A_{n,s}$ and is embedded in $L_{n,s}$ as a deformation retract.

We will change to the language of sphere systems in what follows,  because sphere systems 
detect structure in $M_{n,s}$ which is invisible in terms of graphs.  Let $M=M_{n,s}$, with
boundary spheres $\bdry_0,\cdots,\bdry_{s-1}$. Given a sphere system  $S$
consisting of spheres $S_1,\cdots,S_k$ in
$M$,  the boundary of $M-S$ consists of $\bdry M$ together with pairs $S_i^+$,
$S_i^-$ coming from the $S_i$.  The space $L_{n,s}$ is identified with
systems $S$ such that
 \begin{itemize}
 \item $S$ is simple
\item all spheres in $S$ are nonseparating
 \item the component of $M-S$ which contains  $\bdry_0$
 also contains at least one pair $S_i^+$, $S_i^-$.
\end{itemize}

\section{Sphere systems with connected complement}

We will be interested in sphere systems $S$ in $M = M_{n,s}$ whose complement $M-S$ is
connected. For brevity we call such systems {\it coconnected}. Since a subsystem of a
coconnected system is coconnected, the (isotopy classes of) coconnected systems form a
subcomplex
$Y=Y_{n,s}$ of $S(M)$. 

\begin{figure}[ht]
\centerline
{\relabelbox\small\epsfxsize 3.5truein
\epsfbox{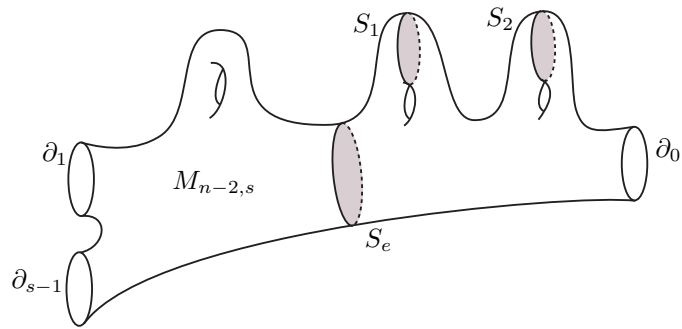}
\relabel{b1}{$\partial_1$}
\relabel{bs}{$\partial_{s-1}$}
\relabel{b0}{$\partial_0$}
\relabel{S1}{$S_1$}
\relabel{S2}{$S_2$}
\relabel{Se}{$S_e$}
\relabel{Mns}{${M}_{n-2,s}$}
\endrelabelbox}
\caption{Enveloping sphere}\label{env}
\end{figure}

When $s >0$ we will also need a refined version of $Y$ whose definition will use the
notion of an {\it enveloping sphere} for a coconnected system $S$. By this we mean a
sphere $S_e$ disjoint from
$S$ which separates $M-S$ into two pieces, one of which is a punctured ball whose
boundary consists of 
$S_e$, the base boundary sphere $\bdry_0$, and the pairs
$S_i^+, S_i^-$ for the spheres
$S_i$ of $S$ (Figure ~\ref{env}). 

When $S$ has $n$ spheres  (the maximum number possible), the
enveloping sphere is unique up to isotopy. If $S$ has $n$ spheres and $s=1$  or
$s=2$, the enveloping sphere will be trivial, bounding a ball if $s=1$ or parallel to
$\bdry_1$ if $s=2$, so in these cases the enveloping sphere is not very interesting.
However, when
$S$ has fewer than $n$ spheres, the enveloping sphere is not unique, and so represents
additional data not specified by $S$.

Let  $V=V_{n,s}$ be the complex
whose simplices are isotopy classes of pairs $(S,S_e)$ where $S$  is a coconnected
system and
$S_e$ is an enveloping sphere for $S$. To obtain a face of the simplex $(S,S_e)$ we
pass to a subsystem
$S'$ of
$S$, with the unique (up to isotopy) enveloping sphere $S'_e$ disjoint from $S$ and
$S_e$.

\begin{theorem} The complex $Y_{n,s}$ is $(n-2)$--connected for $s\geq 0$, and $V_{n,s}$ is
$(n-3)/2$--connected for
$s>0$.
\end{theorem}

\begin{proof}   For $Y_{n,s}$ this is part of Proposition 3.1 of \cite{H}. For $V_{n,s}$ we
will describe a map
$f\colon  SL_{n,s}\to V_{n,s}$ and show this is a homotopy equivalence.

For a sphere system $S$ in $SL_{n,s}$  let $M_0$ denote the component of $M-S$ 
containing $\bdry_0$, and  let $S_1,\cdots,S_k$  be the spheres in $S$ which contribute two
boundary components,
$S_i^+$ and $S_i^-$, to $M_0$. These are the spheres which represent loops
at the basepoint, in the graph picture. The component $M_0$ is a punctured 3--sphere,
as are all components of $M-S$ since $S$ is
 simple.  We define
$f(S)=(T,T_e)$ where $T= \{S_1,\cdots,S_k\} $ and $ T_e$ is the unique enveloping sphere for $T$
which is contained in $ M_0$. 

We claim that $f$ is a poset map, where the poset relation on
$V_{n,s}$ is also  given by inclusion of faces.  Note that any face $S'$ of a sphere system $S$
in $SL_{n,s}$ must contain all of the spheres $S_1,\cdots, S_k$, since removing an $S_i$  results
in a system which is not simple, and therefore not in $SL_{n,s}$.  Since $S'$ is obtained by
removing  spheres other than the $S_i$, the corresponding component $M_0'$  is either the same or
larger than $M_0$.   If $M_0'$ is larger, and $f(S')=(T',T'_e)$, then  the enveloping sphere 
$T'_e$ can be taken to be disjoint from $T_e$, so  $(T',T'_e)$ contains   $(T,T_e)$ as a proper
face. Thus $f$ is indeed a poset map, in which the poset ordering  is reversed.

 Now let $\tau=(T,T_e)$ be an arbitrary simplex
in $V_{n,s}$ and consider the fiber $$f_{\geq \tau}=\{S\in SL_{n,s}\mid f(S)\geq \tau\}.$$ Let $k$ be the number of spheres in $T$. If
$k=n$, then 
$f_{\geq
\tau}=f^{-1}(\tau)$ is a single point (corresponding to a basic thorned rose), since all
spheres in $SL_{n,s}$ are nonseparating.  Now suppose $k<n$.  Then $f_{\geq \tau}$ consists of
all simple systems of nonseparating spheres
 in
$M_{n,s}$ which contain $T$ and are compatible with (disjoint from) $T_e$.  The sphere $T_e$
cuts $M_{n,s}$ into two pieces, one of which is homeomorphic to $M_{n-k,s}$, and $f_{\geq
\tau}$ is isomorphic to the complex of simple systems of nonseparating spheres in this piece.   
This is the space $A_{n-k,s}$, or more precisely the subspace of $A_{n-k,s}$  obtained
by shrinking separating edges. We noted earlier that
$A_{n,s}$  is contractible for all $n>0$, so it follows that the fiber $f_{\geq \tau}$ is
contractible for $k<n$ as well.   Quillen's Poset Lemma \cite{Quillen} then implies that  $f$ is a homotopy
equivalence, so that $V_{n,s}$ is $(n-3)/2$--connected.
\end{proof}

\medskip

The group $ \Gamma_{n,s} $ acts on $Y$ and $ V $, but the stabilizers of simplices are
larger than we would like for the spectral sequence arguments we will be using.  
More specifically, we would prefer that the stabilizers not contain elements that permute the
spheres in a coconnected sphere system or reverse their orientations.  To deal with
orientations we will use oriented versions $Y^{\pm}$ and $V^{\pm}$ of $Y$ and $V$, in which a
$k$--simplex is a
$k$--simplex $S$ or $(S,S_e)$ with a choice of orientation for each sphere of $S$. Thus
a
$k$--simplex of $Y$ or $V$ corresponds to $ 2^{k+1}$ $k$--simplices of $Y^{\pm}$ or
$V^{\pm}$.

\begin{proposition} $Y^{\pm}$ is $(n-2)$--connected and $V^{\pm}_{n,s} $ is $(n-3)/2$--connected.
\end{proposition}

\begin{proof} For $Y^{\pm}$ this was shown in Proposition~3.1 of \cite{H}, and the proof for
$V^{\pm}$ is essentially the same. First choose arbitrarily a positive orientation for each isotopy
class of split-off handles. The subcomplex $V^+$ of $V^{\pm}$ consisting of positively oriented
systems is then a copy of
$V$, so it will suffice to deform a given simplicial map of a triangulated $i$--sphere $S^i$ into  $V^{\pm}$ to a map
with image in $ V^+$, for $ i \le (n-3)/2$. Let $\tau$ be a totally negative simplex of $S^i$ of maximum dimension $k$, where  {\it totally negative} means   $f(\tau)=(T,T_e)$ and all of the spheres in $T$ are negatively oriented. The
linking sphere $S^{i-k-1}$ is either empty (if $k=i$) or  maps to $V^+$, in which case $f$ restricted to $S^{i-k-1}$   can be viewed as a map to
$V^+_{m,s}$ where $M_{m,s}$ is the piece split off by $T_e$.   We have $ m \ge n-k-1$ since $T$ contains at most
$k+1$ spheres. Thus $V^+_{m,s}$ is $(n-k-4)/2$--connected, being
a copy of $V_{m,s}$. The assumption $i\le(n-3)/2$ implies that $i-k-1\le (n-k-4)/2$, so
$f|S^{n-k-1}$ extends to a map $g\co D^{n-k}\to V^+_{m,s}$. (If $k=i$, $g$ just picks a vertex of $V_{m,s}$.)  Writing the star of $\tau$ in $ S^i$
in terms of joins as $
\tau * S^{n-k-1} = \bdry\tau * D^{n-k}$ we can therefore modify $f$ on this star by replacing
it by the map which is the join of $f|\bdry\tau$ with $g$. The modified $f$ is homotopic to
the original
$f$ via the join of $f|\tau$ with $g$. The new $f$ is an improvement on the old one in that no
simplices in the star of $\tau$ of dimension $k$ or greater are now totally negative. Outside the star the map $f$ is
unchanged, so after a finite number of such modifications we will have deformed $f$ to have
image in $V^+$, as desired.
\end{proof}

To eliminate permutations from the simplex stabilizers for the actions of $\G_{n,s}$ on $Y^{\pm}$
and $V^{\pm}$ we make use of a general construction that can be applied to any simplicial complex
$Z$. This construction associates to $Z$ the complex $\Delta(Z)$  whose
$k$--simplices are  the simplicial maps $\Delta^k \to Z$. A $k$--simplex of $\Delta(Z)$ is then an
ordered $(k+1)$--tuple of vertices of $Z$, all of which are vertices of a single simplex
of $Z$. Repeated vertices are allowed in such a $(k+1)$--tuple, so
$\Delta(Z)$ is infinite-dimensional if $Z$ is nonempty. The complex $\Delta(Z)$ is no longer a
simplicial complex but just a
$\Delta$--complex, a CW complex with simplicial cells. There is a canonical map $ \Delta(Z) \to Z $,
and it is a classical fact in algebraic topology this this map induces an isomorphism on homology.
This is shown in Theorem~4.6.8 in \cite{Spanier} for example. [Here is a sketch of the proof. It is
easy to construct a contraction of $\Delta(\Delta^k)$, so the result holds for $\Delta^k$. The result
for finite
$Z$ then follows by induction on the number of simplices via the five lemma and the
Mayer--Vietoris sequence for simplicial homology. A direct limit argument with simplicial
homology then gives the general case.]

It can be shown that the map $\Delta(Z) \to Z$ is not only a homology equivalence, but in fact a homotopy equivalence; however, we will not
need this stronger statement.

We now set $X_{n,s} = \Delta(Y^{\pm}_{n,s})$ and $W_{n,s} = \Delta(V^{\pm}_{n,s})$.

\section{The first stabilization}

Using the spectral sequence associated to the action of the group
$\G=\G_{n+1,s}$ on the complex
$W = W_{n+1,s}$ we will prove:

 \begin{theorem} \label{alpha}
 For $s\geq 1$, the map $\alpha_*\colon H_i(\G_{n,s})\to H_i(\G_{n+1,s})$ induced by
the   inclusion $\alpha\colon\G_{n,s}\to \G_{n+1,s}$ is an isomorphism for $n\geq
2i+2$ and a surjection for $n = 2i+1 $. 
\end{theorem}

Note  that for $s=1$, the map $\alpha$ is the standard inclusion 
$Aut(F_{n})\to Aut(F_{n+1})$, so in this case the theorem gives a slight
improvement on Theorem 7.1 of \cite{HV} which had $2n+3$ in place of $2n+2$.

\begin{proof}  Let $C_*$ be the cellular chain complex of $W$
augmented by a $\Bbb Z$ in degree $ -1$, and let $E\G_*$ be a free $\Z(\G)$--resolution
of $\Z$. Then we have a double complex $C_*\otimes E\G_*$ and hence two spectral sequences associated
to the horizontal and vertical filtrations of this double complex, both converging to the same thing. The spectral sequence obtained by first taking homology with respect to the vertical
boundary maps, in the $E\G_*$ factor, has
$$
E^1_{p,q} = H_q(\G;C_p) = \bigoplus_{\tau}H_q(\G_{\tau};\Z),
$$
where the sum is over all orbits of
$p$--cells $\tau$, $\G_{\tau}$ is the stabilizer of 
$\tau$,   and the second equality  comes from Shapiro's lemma.   Since all elements of $\G_{\tau}$ fix $\tau$ pointwise, the coefficient module $\Z$
is trivial as a $\G$--module, and we will omit this $\Z$ from the notation from now on. 
 If the image of the $p$--simplex $\tau \colon\Delta^p\to V^{\pm}_{n+1,s}$ is the pair $(S,S_e)$ where
the system $S$ has
$k$ distinct spheres, then there is a natural surjection
$\G_{n+1-k,s} \to \G_{\tau}$ coming from the inclusion $M_{n+1-k,s} \to M_{n+1,s}$, where $M_{n+1-k,s}$ is  the part of
$M_{n+1,s}$ on the side of $S_e$ away from $S$.  The surjection $\G_{n+1-k,s} \to \G_{\tau}$ is in
fact an isomorphism by the argument for Lemma~3.3 of
\cite{H}. 
 
The fact that $ W_{n+1,s}$ is $ m $--connected for $m=(n-2)/2$ implies that $ E^{\infty}_{p,q} = 0 $ 
for $ p+q \le m$. This follows by looking at the spectral sequence coming from the horizontal filtration, which first takes 
homology with respect to the boundary maps in $C_*$. This spectral sequence has $E^1_{p,q} = 0 $
for $p \le m$, so since both spectral sequences converge to the same thing, the assertion follows.

We prove the theorem by induction on the homology dimension $i$.  Thus we assume that for $ k < i $, the map $\alpha$ induces an
isomorphism on $H_k$ for $ n \ge 2k+2$ and a surjection for $ n = 2k+1$.

 Let us show first that
$\alpha$ induces a surjection on $H_i$ for $ n \ge 2i+1$. The stabilization $H_i(\G_n) \to
H_i(\G_{n+1})$ is the differential $d^1\colon  E^1_{0,i} \to E^1_{-1,i}$, so to show this is
surjective it suffices to show:
\begin{enumerate}
\item[(i)] $E^{\infty}_{-1,i} = 0 $.
\item [(ii)] $E^2_{p,i-p} = 0 $  for $p>0$. 
  \end{enumerate}
Condition (ii) guarantees that all differentials $d^r$ with $ r > 1 $
going to $E^r_{-1,i}$ must be zero, so only the $d^1$ differential can kill $E^1_{-1,i}$.

It is easy to see that statement (i)  holds since we have already observed that $E^{\infty}_{p,q} = 0 $ for
$ p+q
\le (n-2)/2$, so in particular this holds with $ (p,q) = (-1,i) $ if $ i-1 \le (n-2)/2 $, or $ n \ge 2i$.
To prove (ii) we use the description of the $d^1$ differentials in Chapter 7 of \cite{Brown}. In the
$q$-th row these differentials compute the homology of the quotient $W/\G$ with coefficients in
the system of groups $H_q(\G_\tau)$. For $p+q \le i$ and $ q < i $, and for $\tau'$ a face of $\tau$, the coefficient maps $H_q(\G_\tau)\to H_q(\G_{\tau'})$ induced by inclusion are all isomorphisms by the
induction hypothesis. In addition, induction gives that the coefficient maps from the $ p+q = i+1$ terms to the $ p+q = i $
terms are surjective for
$ q < i $. Thus for $ p+q \le i $ the groups $E^2_{p,q}$ are homology groups of $W/\G$
with locally constant coefficients. In fact, we claim that the coefficients are constant.  
The twisting of
coefficients around the faces of a simplex $\tau$ is induced by diffeomorphisms of $M$ that permute the
nonseparating spheres corresponding to the vertices of $\tau$.   These diffeomorphisms can be chosen
to have support in the part of $M$ cut off by the enveloping sphere for $\tau$.  Since the stabilizer
$\G_\tau$ is supported on the other part of $M$, these diffeomorphisms act trivially
on $H_q(\G_\tau)$.

Next we need the fact that the quotient $W/\G$ is $(n-1)$--connected. This was shown in
Lemma~3.5 of \cite{H} for the analogous quotient $X/\G$, and the two quotients are
combinatorially equivalent. This is because $\G$ acts transitively on $k$--simplices of both the
complexes $V^{\pm}$ and $Y^{\pm}$ from which $W$ and $X$ are constructed, and 
all permutations of the vertices of a simplex of $V^{\pm}$ or $Y^{\pm}$ are realized by the
action of $\G$.

 Thus if $ i \le n-1$ we
conclude that (ii) holds. The condition $ i \le n-1 $ is implied by $ n \ge 2i+1 $, so this
finishes the induction step for surjectivity.

To show that $d^1\colon  E^1_{0,i} \to E^1_{-1,i}$ is injective,  it is enough to show:

\begin{enumerate}
 \item [(iii)] $E^{\infty}_{0,i} = 0 $.

\item [(iv)] $E^2_{p+1,i-p} = 0 $  for $p>0$. 

\item [(v)] The differential $d^1\colon E^1_{1,i} \to E^1_{0,i} $ is zero.
 \end{enumerate}
The arguments for (iii) and (iv) are the same as for (i) and (ii), with a shift of one unit to the right,
corresponding to replacing $ n \ge 2i+1 $ by $n \ge 2i+2$. The connectivity condition for $ W/\G$  is
now
$ i+1
\le n-1 $, and this is implied by $ n \ge 2i+2 $.

For (v) there are two orbits of $1$--cells, one nondegenerate and the other degenerate. The
coefficients are automatically untwisted for the degenerate $1$--cell, and they are untwisted for the
nondegenerate one by the argument given above. In
$W/\G$ each
$1$--cell attaches at each end to the unique
$0$--cell, so it follows that the differential is zero.  
\end{proof}
 
If one examines the proof, one sees that in the inequalities $ n \ge 2i+2$ and $ n \ge 2i+1 $, the
coefficient $2$ in the $2i$ term is forced by the induction process since the $d^2$ differentials go
two units to the left and one unit upward. On the other hand, the constant term $2$ or $1$ comes from
the connectivity $(n-2)/2$ of $W_{n+1,s}$. If this connectivity estimate were
improved, this could yield a corresponding improvement in the constant term of the stable range.

\section{The second stabilization}

We turn now to the
$\beta$ stabilization $H_i(\G_{n,s+2}) \to H_i(\G_{n+1,s}) $. As an easy application of the theorem
in the preceding section we will show:

\begin{theorem}\label{beta} 
 For $s\geq 1$, the map $\beta_*\colon H_i(\Gamma_{n,s+2})\to H_i(\Gamma_{n+1,s})$ induced by  
 $\beta$ is an isomorphism for $n\geq 2i+2$ and a surjection for $ n = 2i+1$.    
\end{theorem}

 \begin{proof}
Surjectivity of $\beta_*$ on $H_i$ follows from surjectivity of $\alpha$ because of
the commutative diagram from Section 2:
$$
\begin{CD}
 \G_{n,s} @>{\alpha}>> \G_{n+1,s}\\
@VV{\mu}V @AA{\beta}A \\
\G_{n,s+1} @>{\mu}>> \G_{n,s+2}
\end{CD}
$$
For injectivity we will first prove the result in an unspecified range $ n \gg  i $ and then
improve this to $ n \ge 2i+2$. We consider the relative homology groups $ H_i$ for the pairs in the
following commutative diagram:
$$
\begin{CD}
(\G_{n,s},\G_{n-1,s+2})@>{\alpha}>>(\G_{n+1,s},\G_{n,s+2})\\
@VV{\mu}V @AA{\beta}A\\
(\G_{n,s+1},\G_{n-1,s+3})@>{\mu}>>(\G_{n,s+2},\G_{n-1,s+4})
\end{CD}
$$
Since $\alpha$ is an isomorphism on the absolute groups $H_i$ for $n\gg i$, it is also an isomorphism
on the relative groups. Then the diagram shows that $\beta$ is surjective on the relative $H_i$ for
$n\gg i$. This is the map on the left in the following commutative square:
$$
\begin{CD}
H_i(\G_{n,s+2},\G_{n-1,s+4})@>>>H_{i-1}(\G_{n-1,s+4}) \\
@V{ \beta_*}VV@V{\beta_*}VV\\
H_i(\G_{n+1,s},\G_{n,s+2})@>>>H_{i-1}(\G_{n,s+2}) \\
\end{CD}
$$
Assuming $ n \gg  i $ always, we can see that the lower map is injective by looking at a portion of the long exact sequence for the pair $(\G_{n+1,s},\G_{n,s+2})$:
$$
H_i(\G_{n,s+2})\  {\buildrel{\beta_*}\over\longrightarrow}\  H_i(\G_{n+1,s}) \longrightarrow H_i(\G_{n+1,s},\G_{n,s+2}) \longrightarrow H_{i-1}(\G_{n,s+2})
$$
The first map is surjective as we saw at the beginning of the proof, so the second map is zero and the third map is injective.   Returning to the commutative square, the composition from upper left to lower right passing through the upper right term is zero since it is two maps in a long exact sequence. Thus the composition through the lower left group is zero. This composition is a surjection followed by an injection, so it follows that $H_i(\G_{n+1,s},\G_{n,s+2})$ is zero. Hence $\beta$ induces an isomorphism on $H_i$  for $ n \gg  i$.

It remains to improve the stable range to what is stated in the theorem. From the first diagram in the
proof, the fact that $\alpha_*$ and
$\beta_*$ are isomorphisms for
$n\gg i$ then
implies that the maps
$\mu_*$  are also isomorphisms for
$n\gg i$, using the fact that $ \mu_* $ is always injective (since $\gamma\circ\mu$ is the
identity, where $\gamma$ is induced by filling in a boundary sphere with a ball). For each 
$s\geq 1$ and each
$i$, we now have a diagram of maps:
$$
\begin{CD}
H_i(\G_{n,s})@>{\alpha_*}>>H_i(\G_{n+1,s})@>{\alpha_*}>>
\cdots@>{\alpha_*}>>H_i(\G_{n+k,s})\\
@VV{\mu_*}V @VV{\mu_*}V @. @VV{\mu_*}V\\
H_i(\G_{n,s+1})@>{\alpha_*}>> H_i(\G_{n+1,s+1})@>{\alpha_*}>>
\cdots@>{\alpha_*}>>H_i(\G_{n+k,s+1})
\end{CD}
$$
In both rows, the horizontal maps are all isomorphisms for $n\geq 2i+2$.  For $k$ large the map 
$\mu_*$ on the far right is an isomorphism.  Since the diagram commutes, we can work
backward to conclude that in fact all the maps
$\mu_*$ in the diagram are isomorphisms, if $n\geq 2i+2$.  
Now the fact that $\beta\circ\mu^2=\alpha$ allows us to conclude  that $\beta_*\colon
H_i(\G_{n,s+2})\to H_i(\G_{n+1,s})$ is an isomorphism for $n\geq 2i+2$, proving the theorem.
\end{proof}

\begin{theorem} The map $\beta_*\colon H_i(\Gamma_{n,2})\to H_i(\Gamma_{n+1,0})$ induced by  
 $\beta$ is an isomorphism for $n\geq 2i+3$ and a surjection for $ n = 2i+2$.    
\end{theorem}

\begin{proof}
This is very similar to the proof of Theorem~\ref{alpha}, using now the action of $ \G
=\G_{n+1,0} $ on the complex
$X=X_{n+1,0}$. The complex $X$ is $(n-1)$--connected,  which is better than the connectivity $
(n-2)/2$ for $W$. The quotient $X/\G$ was shown to be $(n-1)$--connected in \cite{H}. The only
remaining issue to deal with in carrying the proof of Theorem~\ref{alpha} over to the present
context is the justification for why the local
coefficient system on $X/\G$ is untwisted in steps (ii), (iv), and (v) of the proof of Theorem~\ref{alpha}. In the case of
$W/\G$, we used the enveloping spheres of simplices of $W$ for this, but simplices
of $X$ do not have canonical enveloping spheres. 

We will solve this problem by choosing somewhat arbitrary enveloping spheres, and showing
that the map on homology does not depend on this choice, if $n$ is sufficiently
large.   If $\tau$ is a simplex with $k$ distinct spheres, choose an enveloping sphere
which cuts $M_{n+1,0}$ into two pieces, one containing all the spheres in $\tau$, and
one diffeomorphic to $M_{n+1-k,1}$. The coefficients over $\tau$  will be untwisted if
the iterated
$\mu$ map $H_q(\G_{n+1-k,1}) \to H_q(\G_{n+1-k,2k})$ is surjective,  so that $H_q(\G_{\tau})$
is supported on the side of the chosen enveloping sphere away from the spheres of
$\tau$. For step (ii) the worst case we have to deal with is the differential
$E^1_{2,i-1} \to E^1_{1,i-1}$ with
$ H_{i-1}(\G_{n-2,6}) $ in its domain. Here the iterated $\mu$ will be an isomorphism if $ n-2
\ge 2(i-1)+2$, or $ n \ge 2i+2$. Since steps (i) and (ii) give surjectivity, we conclude that
$\beta_*$ is surjective for $n\geq 2i+2$.  

To get injectivity, we need steps (iii)--(v).  For step (iv) we repeat the argument for (ii),
but  shifted one unit to the right, making the inequality $ n \ge 2i+3$. For (v) we are looking
at the differential $E^1_{1,i} \to E^1_{0,i}$ with
$H_i(\G_{n-1,4})$ in its domain, so we need $ n-1 \ge 2i+2$, or $ n \ge 2i+3$. 
\end{proof}

\begin{corollary} The quotient map $Aut(F_n)\to Out(F_n)$ is an isomorphism on homology  in dimension
$i$ for $n\geq 2i+4$, and a surjection for $ n = 2i+3$.
\end{corollary}
\begin{proof} The quotient map $Aut(F_n)
\to Out(F_n)$ is the map $\gamma$ in the following diagram, induced by filling in the boundary sphere with a ball:
$$
\begin{CD}
\Gamma_{n-1,1}@>{\alpha}>>\Gamma_{n,1}\\
@VV{\mu}V @VV{\gamma}V\\
\Gamma_{n-1,2}@>{\beta}>>\Gamma_{n,0}
\end{CD}
$$
 The map $\alpha$ induces an isomorphism on $H_i$ for $n\geq 2i+3$ by Theorem~\ref{alpha}.  It was shown in the proof of Theorem~\ref{beta} that  $\mu$ induces  an isomorphism 
on $H_i$ in the same range.   The corollary follows since $\beta$ induces   an isomorphism for $n\geq 2i+4$ and a surjection for $n=2i+3$. \end{proof}

Because the homology of $Aut(F_n)=\G_{n,1}$ stabilizes, we obtain the following corollary.

\begin{corollary} The group $H_i(Out(F_n))$ is independent of $n$ when $n\geq 2i+4$.
\end{corollary}

\Addresses\recd

   \end{document}